\documentclass[12pt]{amsart}

\usepackage{amssymb,latexsym,amscd}
\usepackage{times}

\newtheorem{theorem}{Theorem}[subsection]
\newtheorem{corollary}[theorem]{Corollary}
\newtheorem{lemma}[theorem]{Lemma}
\theoremstyle{definition}

\theoremstyle{remark}
\newtheorem{remark}[theorem]{Remark}
\newtheorem{claim}[theorem]{Claim}
\numberwithin{equation}{section}

\def\deg{\operatorname{deg}}
\def\Diag{\operatorname{Diag}}
\def\sort{\operatorname{sort}}
\def\trace{\operatorname{trace}}
\def\max_degree{\operatorname{max\_deg}}

\title{A majorization bound for the eigenvalues 
       of some graph Laplacians}

\author{Tamon Stephen}

\address{Institute for Mathematics and its Applications, 
University of Minnesota 
}
\curraddr{Department of Computing and Software,
1280 Main St.~West,
McMaster University,
Hamilton, Ontario,
Canada
L8S 4K1
}
\email{tamon@optlab.mcmaster.ca}

\begin{document}
\maketitle

\begin{abstract}
It is conjectured that the Laplacian spectrum of a graph is majorized
by its conjugate degree sequence.  In this paper, we prove
that this majorization holds for a class of graphs including
trees.  We also show that a generalization of this conjecture to
graphs with Dirichlet boundary conditions is equivalent to the 
original conjecture.
\end{abstract}

% Setup for bibtex
% \nocite{*}  % bib entries to include that are not actually cited.
\bibliographystyle{plain}
% End setup for bibtex.

\section{Introduction}
One way to extract information about the structure of a graph
is to encode the graph in a matrix and study the invariants of 
that matrix, such as the spectrum.  
In this note, we study the spectrum of the ``Combinatorial
Laplacian'' matrix of a graph.

The {\it Combinatorial Laplacian} of a simple graph $G=(V,E)$ on 
the set of $n$ vertices
$V=\{v_1, \ldots, v_n\}$ is the $n \times n$ matrix $L(G)$ defined by:
$$L(G)_{ij} = \begin{cases} \deg(v_i) &\text{if \ } i=j \\ 
                 -1 &\text{if \ } \{i,j\} \in E \\
                 0 &\text{otherwise}
\end{cases}
$$

Here $\deg(v)$ is the {\it degree} of $v$, that is number of edges on $v$.
The matrix $L(G)$ is positive semidefinite, so its eigenvalues are real
and non-negative.
We list them in non-increasing order and with multiplicity:
$$\lambda_1(L(G)) \ge \lambda_2(L(G)) \ge \ldots \ge \lambda_{n-1}(L(G)) \ge \lambda_n(L(G)) = 0$$
When the context is clear, we can write $\lambda_i(G)$ or simply $\lambda_i$.
We abbreviate the sequence of $n$ eigenvalues as $\lambda(L(G))$.

We are interested in the conjecture of Grone and Merris (``GM'')
that the spectrum $\lambda(L(G))$
is majorized by the conjugate partition of the (non-increasing) 
sequence of vertex degrees of $G$ \cite{GM94}.
This question is currently being studied (see for example \cite{DR02}),
but has yet to be resolved.  In this paper, we extend the class of
graphs for which the conjecture is known to hold. 
We also show that if GM holds for graph Laplacians, it also holds
for more general ``Dirichlet Laplacians'' (cf. \cite{CL96})
as conjectured by Duval \cite{Duv03}.

\section{Background and definitions}

\subsection{Graphs}\label{ss:graphs}
Given a graph $G=(V,E)$ with $n=|V|$ vertices and $m=|E|$ edges,
there are several ways to represent $G$ as a matrix.  
There is the {\it edge--incidence
matrix}, a $n \times m$ matrix that records in each column the two 
vertices incident on a given edge.  
For directed graphs we can consider a signed edge--incidence matrix:
$$\partial(G)_{ve} = \begin{cases}
    1 &\text{if \ } v \text{ is the head of edge } e \\
    -1 &\text{if \ } v \text{ is the tail of edge } e \\
    0 &\text{otherwise}
\end{cases}
$$

There is also a $n \times n$ matrix $A(G)$ called the
{\it adjacency matrix} which is defined by:
$$A(G)_{ij} = \begin{cases} 
    1 &\text{if \ } (i,j) \in E \\
    0 &\text{otherwise}
\end{cases}
$$
The diagonal of $A(G)$ is zero.  

We can encode the (vertex) degree sequence of $G$ in non-increasing
order as a vector $d(G)$ of length $n$, and in an $n \times n$ matrix 
$D(G)$ whose diagonal is $d(G)$ and whose off-diagonal elements are 0.
Then the Combinatorial Laplacian $L(G)$ that we study in this
paper is simply $D(G)-A(G)$.  
It is easy to check that if we (arbitrarily) orient $G$ and 
consider the matrix $\partial(G)$ above, we also have 
$L(G)=\partial(G) \partial(G)^t$.

When the graph in question is clear from context, we may
abbreviate the above terms: $L, A, d, D$.

\begin{remark}\label{re:normalized} 
The Laplacian is sometimes defined with entries 
{\it normalized} by dividing by the square roots of the degrees.
However, we do not do that here.
\end{remark}

\subsection{Graph spectra}
The field of spectral graph theory is the study of the
structure of graphs through the spectra (eigenvalues) of
matrices encoding $G$.  Several surveys are available,
including \cite{Chu97} and \cite{CDD80}. 
Besides theoretical aspects of spectral graph theory,
these books describe a wide range of applications of the subject
to chemistry and physics as well as to problems in other 
branches of mathematics such as random walks and isoperimetric problems.

In the case of $L(G)$, there has been considerable effort to study 
the eigenvalue $\lambda_{n-1}$, which is
known as the {\it algebraic connectivity} of $G$.  It can be shown
that $\lambda_{n-1}(L(G)) = 0$ if and only if $G$ is disconnected.
Bounds on $\lambda_{n-1}(L(G))$ then give information on how
well connected a graph is, and are useful, for example, in showing
the existence of expander graphs.  This and other applications
are discussed in \cite{Chu97}.

Currently, little is known about the middle terms of the spectrum.
This is partly because it varies widely depending on the graph.
However, Grone and Merris \cite{GM94} conjecture that the
conjugate partition of the degree sequence majorizes the
spectrum, and showed that the majorization inequalities 
are tight on the class of {\it threshold} graphs.
This conjecture has been extended to simplicial complexes in
recent work by Duval and Reiner \cite{DR02}.

\subsection{Majorization}
We recall that a {\it partition} $p=p(i)$ is a non-increasing 
sequence of natural numbers, and its {\it conjugate} is
the sequence $p^T(j) := |\{i: p(i) \le j\}|$.  
Then $p^T$ has exactly $p(1)$ non-zero elements.
When convenient, we can add or drop trailing zeros in a partition.
For non-increasing real sequences $s$ and $t$ of length $n$, we say
that $s$ is {\it majorized} by $t$ (denoted $s \unlhd t$)
if for all $k \le n$:

\begin{equation}\label{eq:kst}
\sum_{i=1}^{k} s_i \le \sum_{i=1}^{k} t_i
\end{equation}
and 
\begin{equation}\label{eq:nst}
\sum_{i=1}^{n} s_i = \sum_{i=1}^{n} t_i
\end{equation}

The concept of majorization extends to vectors by comparing
the non-increasing vectors produced by sorting the
elements of the vector into non-increasing order.
Given a vector $v$, call the sorted vector $v'$ which contains 
the elements of $v$ sorted in non-increasing order (with multiplicity)
$\sort(v)$. 

In the context of majorization of unsorted vectors, we will
often want to refer to the {\it concatenation} of two vectors
$x$ and $y$ (ie.~the vector which contains the elements of $x$
followed the elements of $y$).  This is denoted $x,y$  as
for example in Lemma~\ref{l:xyz} below.

There is a rich theory of majorization inequalities which 
occur throughout mathematics, see for example \cite{MO81}.  
Matrices are an important source of majorization inequalities.
Notably, the relationship between the diagonal and
spectrum of a Hermitian matrix is characterized by
majorization (see for example \cite{HJ90}).  

We will use the following lemmas about majorization
which can be found in \cite{MO81}:

\begin{lemma}\label{l:xyP}
If $x$ and $y$ are vectors and $P$ is a doubly-stochastic matrix and $x=Py$, 
then $x \unlhd y$.
\end{lemma}

This yields two simple corollaries:
\begin{lemma}\label{l:xyz}
For any vectors $x \unlhd y$ and any vector $z$ we have: $x,z \unlhd y,z$.
\end{lemma}

\begin{lemma}\label{l:ij}
If $x$ and $y$ non-increasing sequences, and $x=y$ 
except that at indices $i<j$ we have $x_i = y_i-a$
and $x_j=y_j+a$ where $a \ge 0$ then $x \unlhd y$.
\end{lemma}
Lemma~\ref{l:ij} says that for
non-increasing sequences transferring units from lower
to higher indices reduces the vector in the majorization 
partial order.  
In particular,
if $x,x',y,y'$ are all {\it non-increasing} sequences, 
$x' \unlhd x$ and $y' \unlhd y$, then
\begin{equation}\label{eq:xyxy}
x' + y'  \unlhd x' + y \unlhd x+y
\end{equation}

Let $A$ and $B$ be positive semidefinite (more generally, Hermitian) 
matrices. Then:
\begin{equation}\label{eq:m1}
\lambda(A),\lambda(B) \unlhd \lambda(A+B)
\end{equation}

A theorem of Fan (1949) says that for positive semidefinite
(more generally, Hermitian) matrices $A$ and $B$:
\begin{equation}\label{eq:m2}
\lambda(A+B) \unlhd \lambda(A)+\lambda(B)
\end{equation}

Let $A$ be an $m \times n$ 0-1 (or incidence) matrix, with row sums
$r_1, \ldots, r_m$ and columns sums $c_1, \ldots, c_n$ both indexed
in non-increasing order.  Let $r^T$ be the conjugate of the 
partition $(r_1, \ldots, r_m)$, and $c$ be the partition
$(c_1, \ldots, c_n)$.  Then the Gale-Ryser theorem asserts that
\begin{equation}\label{eq:gr}
c \unlhd r^T
\end{equation}

\subsection{The Grone-Merris Conjecture}
The Grone-Merris conjecture (GM) is that the spectrum of
the combinatorial Laplacian of a graph is majorized by its
conjugate degree sequence, that is
\begin{equation}\label{eq:gm}
\lambda(G) \unlhd d^T(G)
\end{equation}
Note that 
$$ \sum_{i=1}^{n} d_i^T = \sum_{i=1}^{n} d_i = \trace(L(G))
   = \sum_{i=1}^{n} \lambda_i $$

If we ignore isolated vertices (which contribute only zero
entries to $\lambda$ and $d$) we will have $d^T_1=n$.
Using this fact, it is possible to show that
\begin{equation}\label{eq:1ld}
\lambda_1 \le d^T_1
\end{equation}
Three short proofs of this are given in \cite{DR02}.
The authors then continue to prove the second majorization inequality
\begin{equation}\label{eq:2ld}
\lambda_1 + \lambda_2 \le d^T_1 + d^T_2
\end{equation}
However, their proof would be difficult to extend.

There are several other facts which fit well with the GM
conjecture.  One is that if the GM conjecture holds, then 
the instances where (\ref{eq:gm}) holds with equality are well-understood,
these would be the threshold graphs of Section~\ref{ss:th}.
Also, since $d$ and $\lambda$ are respectively the diagonal and spectrum
of $L(G)$ we have $d \unlhd \lambda$.  Combining this with GM 
gives $d \unlhd d^T$, a fact that has been proved combinatorially.
We refer to \cite{DR02} for further discussion.

\begin{remark}[Complements]\label{re:complements}
Given a graph $G$, we can study its {\it complement} $\overline{G}$,
the graph whose edges are exactly those not included in $G$.
For a graph $G$ with $n$ vertices the $i$th largest vertex of $G$
is the $(n-i)$th largest vertex of $\overline{G}$, and we have
$d_i(G)=n-1-d_{n-i}(\overline{G})$.  Translating this to 
the conjugate partition $d^T$ yields:
$d_i^T(G)=n-d_{n-1-i}^T(\overline{G})$ with
$d_n^T(G)=d_n^T(\overline{G})=0$.

The relationship between $\lambda(G)$ and $\lambda(\overline{G})$
is the same as between $d_n^T(G)$ and $d_n^T(\overline{G})$.
This follows from the fact that
$L(G)+L(\overline{G})=n I_n-J_n$ where $J_n$ is the $n \times n$
matrix of ones.  We observe that the matrix $n I_n-J_n$ sends the special
eigenvector $e_n$ ($n$ ones) to 0, and acts as the scalar $n$
on $e_n^\perp$.  
Both $L(G)$ and $L(\overline{G})$ also send $e_n$ to 0, giving us
$\lambda_n(G)=\lambda_n(\overline{G})=0$.
Since $L(G)$ and $L(\overline{G})$ sum to $n I_n$ on $e_n^\perp$
they have the same set of eigenvectors on $e_n^\perp$, and 
and for each eigenvector the corresponding eigenvalues for 
$L(G)$ and $L(\overline{G})$ sum to $n$.
Thus 
$\lambda_i(G)=n-\lambda_{n-1-i}(\overline{G})$.
As a consequence, GM holds for $G$ if and only if GM holds for $\overline{G}$.
\end{remark}

\section{Grone-Merris on classes of graphs}
In this section we
give further evidence for the Grone-Merris conjecture by remarking that
it holds for several classes of graphs including threshold graphs,
regular graphs and trees.

\subsection{Threshold graphs}\label{ss:th}
The GM conjecture was originally formulated in the context
of {\it threshold} graphs, which are a class of graphs with several
extremal properties.  An introduction to threshold graphs is \cite{MP95}.  
Threshold graphs are the graphs
that can be constructed recursively by adding isolated vertices
and taking graph complements.  It turns out that they are also
characterized by degree sequences: the convex hull of possible 
(unordered) degree sequences of an $n$ vertex graph defines
a polytope.  
The extreme points of this polytope are the degree sequences that have
a unique labelled realization, and these are exactly the threshold graphs.

Threshold graphs are interesting from the point of view of
spectra.  Both Kelmans and Hammer \cite{HK96} and Grone
and Merris \cite{GM94} investigated the question of
which graphs have integer spectra.  They found that threshold
graphs are one class of graphs that have integer spectra and
showed for these graphs that $\lambda(G)=d^T(G)$.  

In the process of showing this equality for threshold graphs,
Grone and Merris observed that for non-threshold graphs, the 
majorization inequality $\lambda(G) \unlhd d^T(G)$
appears to hold, and made their conjecture.  
We could describe the conjecture as saying that threshold 
graphs are extreme in
terms of spectra, and that the these extreme spectra can be
interpreted as conjugate degree sequences.

\subsection{Regular and nearly regular graphs}
For some small classes of graphs, it can be easily shown
that the GM conjecture holds.
Consider a $k$-regular graph $G$ on $n$ vertices
(in a {\it $k$-regular} graph, all vertices have degree $k$).
Then the degree sequence $d(G)$ is $k$ repeated $n$ times, and its
conjugate $d^T(G)$ is $n$ repeated $k$ times followed by $n-k$ zeros.
Thus $d^T$ majorizes every non-negative sequence of sum $kn$ whose largest
terms is at most $n$, and in particular $\lambda \unlhd d^T$.
Indeed, this proof shows that GM holds for what we might call 
{\it nearly regular} graphs, that is graphs whose vertices have
degree either $k$ or $(k-1)$.

\subsection{Graphs with low maximum degree}
Using facts about the initial GM inequalities
we can prove that GM must hold for graphs with
low maximal degree.  For example, if a graph has maximum vertex
degree 2, then $d^T_3 = d^T_4 = \ldots = d^T_n = 0$, so for
$k=2,3,\ldots,n$:
$$\sum_{i=1}^{k} \lambda_i \le 
  \sum_{i=1}^{n} \lambda_i = \sum_{i=1}^{n} d^T_i = \sum_{i=1}^{k} d^T_i$$
More generally, the GM inequalities for $k \ge \max_degree(G)$ hold
trivially.  Thus GM holds for graphs
of maximum degree 2 by (\ref{eq:1ld}).
Using Duval and Reiner's result (\ref{eq:2ld}), we get that
GM holds for graphs of maximum degree 3.

\subsection{Trees and more}
It is tempting to try to prove GM inductively by breaking 
graphs into simpler components on which GM clearly holds.
In this section, we show that if $G$ is ``almost'' the 
union of two smaller graphs on which GM holds then GM holds
for $G$ as well.  We apply this construction to show that
GM holds for trees.

Take two graphs $A=(V_A,E_A)$ and $B=(V_B,E_B)$ 
on disjoint vertex sets $V_A$ and $V_B$.
Define their {\it disjoint sum} to be 
$A+B=(V_A \cup V_B, E_A \cup E_B)$.  Assuming $V_A$ and $V_B$
are not empty this is a disconnected graph.
Now take two graphs $G=(V,E_G)$ and $H=(V,E_H)$
on the same vertex set $V$.  
Define their {\it union} as $G \cup H = (V, E_G \cup E_H)$.

Given the spectra and conjugate degree sequences of $A$ and $B$,
the spectrum of $A+B$ is (up to ordering) 
$\lambda(A+B) = (\lambda(A), \lambda(B))$,
while the conjugate degree sequence of $A+B$ is $d^T(A+B) = d^T(A)+d^T(B)$
(taking each vector to have length $n$).  Thus by~\ref{eq:m1}
if $\lambda(A) \unlhd d^T(A)$ and $\lambda(B) \unlhd d^T(B)$
we will have $\lambda(A+B) \unlhd d^T(A+B)$.

In a typical situation, where neither $A$ or $B$ is very
small, we would expect the above majorization inequality
to hold with considerable slack.  We can use this slack to
show that if we add a few more edges to $A+B$ the majorization
will still hold.  

\begin{theorem}\label{th:abc}
Take graphs $A$ or $B$ on disjoint vertex sets $V_A$ and $V_B$.
Let $G = A+B$ and on $V = V_A \cup V_B$ let 
$C$ be a graph of ``new edges'' between $V_A$ and $V_B$.
Assume that GM holds on $A$, $B$ and $C$, ie.~that
$\lambda(A) \unlhd d^T(A)$, $\lambda(B) \unlhd d^T(B)$
and $\lambda(C) \unlhd d^T(C)$.
Additionally, assume that $d^T_i(C) \le d^T_i(A), d^T_i(B)$ for all $i$,
and that $d^T_1(B) \le d^T_m(A)$ where $m$ is the largest non-zero
index of $d^T(C)$ (equivalently, $m$ is the maximum vertex degree in $C$).
Let $H=C \cup G$.  Then:
\begin{equation}\label{eq:mabc}
\lambda(H) \unlhd d^T(H)
\end{equation}
\end{theorem}

\begin{proof}
Let $k$ be the larger of $\max_degree(A)$ and $\max_degree(B)$.  Note that 
$$d^T(G) = d^T(A)+d^T(B) = 
(d^T_1(A)+d^T_1(B), d^T_2(A)+d^T_2(B), \ldots, d^T_k(A)+d^T_k(B), 0, \ldots, 0)
$$
\begin{claim}\label{c:gc}
$$ d^T(H) \unrhd (d^T_1(G), d^T_2(G), \ldots, 
          d^T_k(G), d^T_1(C), \ldots, d^T_m(C) )$$
\end{claim}

\noindent
{\it Proof of Claim.}
The term on the right is the concatenation of two partitions,
$d^T(G)$ and $d^T(C)$.  The columns of $d^T(G)$ index the vertices
of $G$ and the length of a column gives the degree of the corresponding
vertex.  Since this claim is purely about the combinatorics of
of degree sequences, we introduce a series of
intermediate ``partial graphs'' where edges are allowed to have only
one end.  Degree sequences and their conjugates are still well defined
for such objects.

Let $G_0=G$ and $C_0=C$.  Define $G_i$ by moving one end of an edge
from every non-isolated vertex of $C_{i-1}$ to $G_{i-1}$, and let 
$C_{i}$ contain whatever is left.  Iterating this, for
some $l \ge 0$ we will have $G_l=H$ and $C_l$ consisting entirely
of isolated vertices.  
Then the claim will follow if we can show that:
$$ d^T(G_0), d^T(C_0) \unlhd d^T(G_1), d^T(C_1) \unlhd \ldots
 \unlhd d^T(G_l), d^T(C_l) $$
Compare the partitions at the $(i-1)$st majorization: we 
remove the first row of $d^T(C_{i-1})$ and put each element from
that row into a separate column (representing a distinct vertex in
$G$) of $d^T(G_{i-1})$.  Where there are columns of equal length
in $d^T(G_{i-1})$ 
they should be ordered so that those acquiring new elements come
first.
To see that this operation increases the partition
in the majorization partial order, observe that
after ignoring the (unchanged) contents of $d^T(C_i)$ it is
equivalent to sorting the new row into the partition,
using Lemma~\ref{l:ij} to move its final (rightmost) element 
to the proper column and repeating as necessary.
\vspace{3mm}

This completes the proof of the Claim~\ref{c:gc}, and gives us:
$$ d^T(H) \unrhd (d^T_1(A)+d^T_1(B), d^T_2(A)+d^T_2(B), \ldots, 
          d^T_k(A)+d^T_k(B), d^T_1(C), \ldots, d^T_m(C) )$$
If we sort the vector on the right into non-increasing
order, the first $m$ terms will remain fixed by the assumptions
that $d^T_m(A) \ge d^T_1(B) \ge d^T_1(C)$.  
Since we have assumed that $d^T_i(C) \le d^T_i(B)$ for all $i$, 
we can apply Lemma~\ref{l:ij} to the reordered sequence to get:
\begin{align*}
d^T(H)  \unrhd 
 (d^T_1(A)+d^T_1(C),& d^T_2(A)+d^T_2(C), \ldots, 
          d^T_m(A)+d^T_m(C), \\ & d^T_{m+1}(A), \ldots, d^T_{k}(A),
	d^T_1(B), \ldots, d^T_k(B) )
\end{align*}
The right hand term decomposes as:
$$(d^T_1(A), \ldots, d^T_k(A), d^T_1(B), \ldots, d^T_k(B))
 + (d^T_1(C), \ldots, d^T_m(C),0,\ldots,0)$$
Since we assume $d^T_m(A) \ge d^T_1(B)$, the first $m$ entries
of $(d^T(A), d^T(B))$ will remain unchanged if the vector is sorted.
Thus:
\begin{equation}\label{eq:dt}
d^T(H) \unrhd \sort(d^T(A),d^T(B)) + d^T(C)
\end{equation}
By (\ref{eq:xyxy}) we can apply the majorizations of $\lambda$
by $d^T$ for $A, B, C$ to the above terms to get:
$$\begin{array}{cl}
d^T(H) & \unrhd ~ \sort(d^T(A),d^T(B)) + d^T(C) 
           ~ \unrhd ~ \sort(d^T(A),d^T(B)) + \lambda(C)\\
       & \unrhd ~ \sort(d^T(A),\lambda(B)) + \lambda(C)
          ~  \unrhd ~  \sort(\lambda(A),\lambda(B)) + \lambda(C)
\end{array}$$
The two terms on the right side of this equation are spectra
of $L(G)$ and $L(C)$ respectively.  Hence by Fan's theorem (\ref{eq:m2})
their sum majorizes the spectrum of L(G)+L(C)=L(H):
$$d^T(H) \unrhd \lambda(G)+ \lambda(C) \unrhd \lambda(H)$$
\end{proof}

More generally, we could replace the conditions in the statement 
of Theorem~\ref{th:abc} with the condition (\ref{eq:dt}), 
which can be checked combinatorially.
The conditions in the theorem statement and equation (\ref{eq:dt})
are most likely to be satisfied if $C$ is small relative to $A$ and $B$.

A useful case is when $C$ consists of $k$ disjoint edges. 
Then $m=1$ and $d^T_1(C)=2k$.  Without loss of generality
we can take $d_1(A) \ge d_1(B)$ and the only condition that we will
need to check is that $d_1(A), d_1(B) \ge d_1(C)$, ie.~both $A$ and
$B$ must have at least $2k$ non-isolated vertices.

The strategy for applying Theorem~\ref{th:abc} to show that a given 
graph $H$ satisfies GM is to find a ``cut'' $C$ for it that
contains few edges and divides $H$ into relatively large components.
For example we have the following result:

\begin{corollary}
The GM conjecture holds for trees.
\end{corollary}

\begin{proof}
Proceed by induction on the diameter of the graph.
If $T$ has diameter 1 or 2, then there is a vertex $v$
which is the neighbour of all the remaining vertices
and $T$ is a threshold graph.  
So GM holds with equality for $T$.

Otherwise, we can find some edge $e$ that does not have a leaf vertex.
Since $T$ is a tree, $e$ is a cut edge and divides $T$ into two
non-trivial connected components, $A$ and $B$.  We apply induction to
$A$ and $B$ and apply Theorem~\ref{th:abc} to $H=(A+B) \cup C$
where $C$ is the graph on the vertex set of $T$ containing the
single edge $e$.
\end{proof}

\begin{remark}[Small Graphs]  
The facts in this section allow us to check that GM holds for
some small graphs without directly computing eigenvalues.
For example, since the GM condition is closed under complement 
(see~\ref{re:complements}) for graphs on up
to 5 vertices it is enough to observe that either $G$ or $\overline{G}$
has maximum degree $\le 3$.  Out of 156 graphs on 6 vertices, 146
can be decomposed into smaller graphs $(A+B) \cup C$ using
Theorem~\ref{th:abc}.  Calculating the eigenvalues of the remaining
10 does not yield a counterexample.
\end{remark}

\section{Simplices and pairs}
The most recent work relating to the GM conjecture has been to study the
spectra of more general structures than graphs, such as simplicial
complexes and simplicial family pairs.
In this section we show that the generalization of GM to graphs with
Dirichlet boundary conditions is equivalent to the original conjecture 
and may be useful in approaching GM.

\subsection{Simplicial complexes}
In \cite{DR02}, the authors look at {\it simplicial complexes},
which are higher dimensional analogues of simple graphs (see for example
\cite{Mun84}).
A set of faces of a given dimension $i$
is called an {\it i-family}.  Given a simplicial complex $\Delta$
we can denote the $i$-family of all faces in $\Delta$ of dimension $i$
as $\Delta^{(i)}$.
For example, a graph is a 1-dimensional complex, and its edge set is
the 1-family $\Delta^{(1)}$.
Define the degree sequence $d$ of an $i$-family to be
the list of the numbers of $i$-faces from the family
incident on each vertex, and sorted into non-increasing order.  
We can then define $d(\Delta, i)$ as the degree sequence of
$\Delta^{(i)}$, which we can abbreviate to $d(\Delta)$ or $d$
when the context is clear.

We define the {\it chain group} $C_i(\Delta)$ of formal linear
combinations of elements of $\Delta^{(i)}$, and generalize 
the signed incidence matrix $\partial$ of Section~\ref{ss:graphs}
to a signed boundary map 
$\partial_i:C_i(\Delta) \rightarrow C_{i-1}(\Delta)$.
This allows us to define a {\it Laplacian} on $C_i(\Delta)$, 
namely $L_i(\Delta) = \partial_i \partial_i^T$,
and study its corresponding spectrum $s(\Delta,i)$ sometimes
abbreviated $s(\Delta)$ or $s$.

Duval and Reiner \cite{DR02} looked at 
{\it shifted} simplicial complexes, which are a
generalization of threshold graphs to complexes.  
They showed that for a shifted complex $\Delta$ and any $i$,
we have $s(\Delta,i)=d^T(\Delta,i)$.  They then conjectured
that GM also holds for complexes, ie.~that for any complex
and any $i$ we have:
\begin{equation}\label{eq:gmc}
s(\Delta,i) \unlhd d^T(\Delta,i)
\end{equation}
They also show that some related facts, such as equation
(\ref{eq:1ld}) generalize to complexes.

\subsection{Simplicial pairs}
In \cite{Duv03}, Duval continues by studying {\it relative (family) pairs} 
$(K,K')$ where the set $K=\Delta^{(i)}$ for some $i$ 
is taken modulo a family of $(i-1)$-faces $K' \subseteq \Delta^{(i-1)}$.
When $K'= \emptyset$, this reduces to the situation of the previous
section.

\begin{remark}\label{re:dirchlet}
In the case $i=1$ this is the edge set of a graph ($K$) with a set of 
{\it deleted} boundary vertices $K'$.  An edge attached to a deleted
vertex will not be removed -- it remains as part of the pair,
but we now think of the edge as having a hole on one (or both) ends.

This type of graph with a boundary appears in conformal invariant
theory.  In this language, the relative Laplacian of an
(edge, vertex) pair is sometimes referred to as a 
{\it Dirichlet Laplacian} and its eigenvalues as 
{\it Dirichlet eigenvalues}, see for example \cite{CL96}.
Recently \cite{CE02} used the spectrum of the Dirichlet Laplacian in 
the analysis of ``chip-firing games'', which are processes on graphs
that have an absorbing (Dirichlet) boundary at some vertices.
\end{remark}

We can form chain groups $C_i(K)$ and $C_{i-1}(K,K')$
and use these to define a (signed) boundary operator on
the pair $\partial(K,K'): C_i(K) \rightarrow C_{i-1}(K,K')$.
Hence we get a Laplacian for family pairs 
$L(K,K')= \partial(K,K') \partial(K,K')^T$.
Considered as a matrix, $L(K,K')$ will be the principal submatrix
of $L(K)$ whose rows are indexed by the $i$-faces in $\Delta^{(i-1)}-K'$.
Finally, we get a spectrum $s(K,K')$ for family pairs from the 
eigenvalues of $L(K,K')$.

Duval defines the degree $d_v(K,K')$ of vertex $v$ 
(in the case of a graph, $v$ is allowed to be in $K'$)
relative to the pair $(K,K')$ as the number of faces in $K$
that contain $v$ such that $K-\{v\}$ is in $\Delta^{(i-1)}-K'$.
This allows him to define the degree sequence $d(K,K')$ for
pairs, and to conjecture that GM holds for relative pairs:
\begin{equation}\label{eq:gmr}
s(K,K') \unlhd d^T(K,K')
\end{equation}

\subsection{The Grone-Merris conjecture for relative pairs}
It turns out that at least in the case of (edge, vertex) pairs
that (\ref{eq:gmr}) follows from the original GM conjecture for
graphs.

\begin{theorem}\label{th:gmgmp}
GM for graphs $\Rightarrow$ GM for (edge, vertex) pairs.
\end{theorem}

\begin{proof}
Let $G=(V,E)$ be a graph with $D \subseteq V$ a set of ``deleted''
vertices.  Let $U=V - D$ be the remaining ``undeleted'' vertices.
We will assume that GM holds only on the undeleted part of the
graph, ie.~$G|_U$.  So we have $s(G|_U) \unlhd d^T(G|_U)$.
We can ignore the edges in $G|_D$ completely, since they have
no effect on either $s(G)$ or $d(G)$.
The remaining edges connect vertices in $D$ to vertices in $U$.
Define $G'$ to be the graph on $V$ whose edge are exactly the
edges of $G$ between $D$ and $U$.  Let $a$ be the degree sequence
of the deleted vertices in $G'$ and $b$ be the degree sequence
of the undeleted vertices in $G'$.

We can compute $d^T(E,D)$ in terms of
the degree sequences and spectra of $G|_U$, $G'$ and $G|_D$
since $d^T_i(E,D)$ is the number of vertices (deleted or not)
attached to at least $i$ non-deleted vertices.
The number of such vertices in $U$ will be $d^T_i(G|_U)$,
and the number in $D$ will be $d^T_i(G')=a^T$.
Hence $d^T(E,D) = d^T_i(G|_U) + a^T$.

Now consider the Laplacian $L(E,D)$.
This is the submatrix of $L(G)$ indexed by $U$.
An edge $(i,j)$ in $G|_U$ contributes to entries
$ii, ij, ji, jj$ in both $L(E,D)$ and $L(G)$.
An edge in $G'$, say from $i \in U$ to $j \in D$ contributes
only to entry $ii$, and an edge in $G|_D$ does not affect
$L(E,D)$.  So, we have $L(E,D) = L(G|_U) + \Diag(b)$, and
by (\ref{eq:m2}) we have:
\begin{equation}\label{eq:speced}
s(E,D) \unlhd s(G|_U) + b
\end{equation}

We complete our equivalence by appealing to the Gale-Ryser 
theorem (\ref{eq:gr}) to claim that $b \unlhd a^T$.
This follows from the fact that $a$ and $b$ are row and
column sums (in non-increasing order) of the $|D| \times |U|$
bipartite incidence matrix for $G'$.  Combining with the
assumption that $s(G|_U) \unlhd d^T(G|_U)$ and 
(\ref{eq:speced}) we get:
$$s(E,D) \unlhd s(G|_U) + b \unlhd d^T(G|_U) + a^T = d^T(E,D)$$
\end{proof}

This proof relies on the bipartite structure of $G'$, so
it is not immediately obvious how to extend it 
to higher dimensional complexes.  It would be interesting to
do this.

\begin{remark}
Because the induction used to prove Theorem~\ref{th:gmgmp}
requires only that the ``undeleted'' part of the graph satisfy
GM, it is tempting to attack the original GM conjecture by
showing if GM holds for a pair $(G,\{v\})$ then GM holds for $G$.
\end{remark}

\section{Acknowledgements}
Thanks to Vic Reiner for introducing me to this question
and for comments throughout.  

\bibliography{gm}

\end{document}